\documentclass[a4paper]{amsart}

\theoremstyle{plain}    
\newtheorem{thm}{Theorem}[section]
\numberwithin{equation}{section} 
\numberwithin{figure}{section} 
\theoremstyle{plain}    
\newtheorem{cor}[thm]{Corollary} 
\theoremstyle{plain}    
\newtheorem{lem}[thm]{Lemma} 
\newtheorem{fact}[thm]{Fact} 
\theoremstyle{plain}    
\newtheorem{prop}[thm]{Proposition} 
\theoremstyle{plain}    
\newtheorem{conjecture}[thm]{Conjecture} 
\theoremstyle{definition}
\newtheorem{problem}[thm]{Problem}
\theoremstyle{remark}
\newtheorem{rem}[thm]{Remark}
\theoremstyle{remark}    
\newtheorem{notation}[thm]{Notation} 
 
\newtheorem{assumption}[thm]{Assumptions} 

\usepackage{amscd}
\usepackage{times}

\usepackage[arrow,matrix,curve,ps]{xy}
\xyoption{dvips}
\CompileMatrices

\renewcommand{\O}{\mbox{$\mathcal{O}$}}
\renewcommand{\P}{\mbox{$\mathbb{P}$}}
\newcommand{\C}{\mbox{$\mathbb{C}$}}
\newcommand{\Z}{\mbox{$\mathbb{Z}$}}
\newcommand{\mc}{\mathcal}

\DeclareMathOperator{\Alb}{Alb}
\DeclareMathOperator{\Aut}{Aut}
\DeclareMathOperator{\codim}{codim}
\DeclareMathOperator{\Div}{Div}
\DeclareMathOperator{\End}{End}
\DeclareMathOperator{\Ext}{Ext}
\DeclareMathOperator{\Hom}{Hom}
\DeclareMathOperator{\locus}{locus}
\DeclareMathOperator{\Pic}{Pic}
\DeclareMathOperator{\rk+}{rank^+}
\DeclareMathOperator{\rank}{rank}
\DeclareMathOperator{\Sing}{Sing}
\DeclareMathOperator{\Supp}{Supp}

\newcommand{\Lbullet}{{{\bf L}^\bullet}}

\begin{document}

\title{Projective Contact Manifolds} 

\author[S.~Kebekus, T.~Peternell, A.~J.~Sommese and
J.~A.~Wi{\'s}niewski] {Stefan~Kebekus, Thomas~Peternell,
  Andrew~J.~Sommese and Jaros{\l}aw~A.~Wi{\'s}niewski}

\address{Stefan Kebekus, Lehrstuhl Mathematik VIII, Universit\"at
  Bayreuth, 95440 Bayreuth, Germany, {\tt
    stefan.kebekus@uni-bayreuth.de} \newline \indent
  {\em current address:} Stefan Kebekus, Research Institute for Mathematical 
  Studies, Kyoto University, Kyoto 606-01, Japan
}

\address{Thomas Peternell, Lehrstuhl Mathematik I, Universit\"at
  Bayreuth, 95440 Bayreuth, Germany, {\tt
    thomas.peternell@uni-bayreuth.de}}

\address{Andrew J.~Sommese, Department of Mathematics, University of
  Notre Dame, Notre Dame, Indiana 46556-5683, USA, {\tt
    sommese@nd.edu}}

\address{Jaros{\l}aw A. Wi{\'s}niewski, Instytut Matematyki UW,
  Banacha 2, 02-097 Warszawa, Poland, {\tt jarekw@mimuw.edu.pl}}

\thanks{At the time of writing the paper, J.Wisniewski was Humboldt-Stipendiat
  at G\"ottingen. He would also like to acknowledge partial support
  from Polish KBN}

\date{\today}

\begin{abstract}
  The present work is concerned with the study of complex projective
  manifolds $X$ which carry an additional complex contact structure.
  In the first part of the paper we show that if $K_X$ is not nef,
  then either $X$ is Fano and $b_2(X)=1$, or $X$ is of the form
  $\P(T_Y)$, where $Y$ is a projective manifold. In the second part of
  the paper we consider contact manifolds where $K_X$ is nef or, more
  generally, complex projective manifolds where the bundle of
  holomorphic 1-forms contains a nef sub-bundle of rank 1.
\end{abstract}


\maketitle
\tableofcontents

\section{Introduction}

A compact complex manifold $X$ of dimension $2n+1$ together with a
subbundle $F\subset T_X$ of rank $2n$ is a \emph{contact manifold} if
the pairing $F\times F\to T_X/F=:L$ induced by the Lie bracket is
everywhere non-degenerate. Alternatively, the non-degeneracy condition
can be reformulated as follows: if $\theta \in H^0(X,\Omega_X^1\otimes
L)$ is the induced form, then an elementary calculation shows that
$\theta\wedge (d\theta)^{\wedge k}$ gives a section of
$\Omega^{2k+1}_X\otimes L^{k}$. The subbundle $F=\ker(\theta)$ is
non-degenerate in the above sense if $\theta \wedge (d\theta )^n$ has
no zeroes.

Contact structures first came up in real geometry. One of the
interests in complex geometry lies in the connection with twistor
spaces and quaternionic K{\"a}hler manifolds.  We refer to
\cite{Bea99} and \cite{Leb95} for excellent introductions to these
matters.

This paper contributes to the classification of projective contact
manifolds. It has been proved by \cite{Dru98} that $\kappa (X) =
-\infty$. Therefore it seems natural to apply Mori theory for the
classification.

In section \ref{morichap:jar} we suppose that $K_X$ is not nef,
i.e.~that there is a curve $C\subset X$ such that $K_X.C<0$. The main
result of this section is:

\begin{thm}\label{mainthm}
  Let $X$ be a projective contact manifold and assume that $K_X$ is
  not nef. Then either $X$ is Fano with $b_2=1$, or there exists a
  smooth projective manifold $Y$ such that $(X,L) \cong
  (\P(T_Y),\O_{\P(T_Y)}(1))$.
\end{thm} 

There are two further problems to be addressed for a complete
classification. These are the following

\begin{problem}
  It is necessary to classify the Fano contact manifolds with
  $b_2(X)=1$.
\end{problem}

It is conjectured that a Fano contact manifold with $b_2=1$ is a
homogeneous variety which is the unique closed orbit in the
projectivized (co)adjoint representation of a simple algebraic Lie
group. The results of Beauville \cite{Bea98} provide a strong
evidence for this conjecture. The present work, however,
is not concerned with this problem as we are primarily interested in
those $X$ which have second Betti-number $>1$.

\begin{problem}
One has to show that $K_X$ cannot be nef. 
\end{problem}
This should be a general consequence of $\kappa (X)=-\infty $; it is
expected that manifolds with negative Kodaira dimension are uniruled
(in dimension 3 this has already been shown by Miyaoka and Mori).
However this conjecture is completely open in dimension at least 4.

In section \ref{KxNEF_CHAP} we prove that a hypothetical projective
contact manifold $X$ with $K_X$ nef has $K^2_X \equiv 0$. We present
an approach that can help to rule out the case that $K_X$ is nef: one
knows that $T_X$ is semi-stable with respect to the (degenerate)
polarization $(K_X,H^{m-2})$, where $m=\dim X$ and $H$ an arbitrary
ample divisor. We would need to have semi-stability with respect to
$(K_X+\epsilon H,H^{m-2})$ for some small $\epsilon >0$.

This problem motivates us to consider a more general situation in
section \ref{nef_in_Omega_CHAP}, namely that $\Omega _X^1$ contains a
nef locally free subsheaf, $L^*\subset \Omega ^1_X$ which is
proportional to $K_X$. In this case we have $\kappa (X)\leq 1$.  If
$\kappa (X)=1$, then $K_X$ is semi-ample and we obtain information on
$L^*$ via the Iitaka fibration. If $\kappa (X)=0$ we conjecture that
$X$ is a quotient of a product $A\times Y$ where $A$ is Abelian and
$Y$ is simply connected. We prove this conjecture if $\dim X\leq 4$ or
if Ueno's Conjecture K holds.

For the benefit of readers coming from differential geometry, we have
added at various places some explanations concerning basic concepts of 
algebraic geometry, e.g.~Mori theory.

For interesting discussions on contact manifolds we would like to
thank A.Beauville and C. Okonek.

\section{Contact manifolds where $K_X$ is not nef}
\label{morichap:jar}

\begin{notation}
  Let $X$ be manifold of dimension $2n+1$ with a line bundle $L$.  In
  this section the tensor product ${\mathcal F}\otimes L^{\otimes k}$
  of an $\O_X$-sheaf ${\mathcal F}$ with a tensor power of $L$ will be
  denoted by ${\mathcal F}(k)$ so that by $\Omega_X^p(k)$ we denote
  the sheaf of holomorphic $p$-forms on $X$ twisted by $L^{\otimes
    k}$.
  
  Let $\theta$ be a twisted holomorphic 1-form $\theta\in
  H^0(X,\Omega^1_X\otimes L)=H^0(X,\Omega^1_X(1))$.
\end{notation}

For a trivializing covering $\{U_i\}$ in which the bundle $L$ is given
by transition functions $g_{ij}$, the form $\theta$ will be locally
trivialized by $\theta_i\in\Omega^1_X(U_i)$ with the relation
$\theta_i=g_{ij}\cdot \theta_j$.  It is immediate to see that although
usually $d\theta_i$ do not glue to a section of $\Omega_X^2(1)$, yet
$\theta_i\wedge(d\theta_i)^n$ define a section of
$\Omega_X^{2n+1}(n+1)$ which we call just $\theta \wedge (d\theta)^n$.

The manifold $X$ is a {\em contact manifold} with the contact form
$\theta\in H^0(X,\Omega^1_X(1))$ if $\theta\wedge (d\theta)^n$ does not
vanish anywhere; this implies that 
$$
L^{\otimes(n+1)}\cong \det T_X = -K_X,
$$
where $K_X$ is the canonical bundle. The contact form $\theta$ defines 
a vector bundle epimorphism
$$
T_X \to L \to 0.
$$
We let $F$ denote its kernel.

For an excellent introduction to contact manifolds we refer the reader
to \cite{Leb95}.

\subsection{Preliminaries, $\C^*$-bundles and the Atiyah extension
  class}
\label{jar:2.1}

First, let us discuss generalities related to a definition of the
Chern class of a line bundle in terms of extensions of sheaves of
differentials. This approach is due to Atiyah \cite{Ati57}; we recall
it and adjust it to our particular set-up. In the subsequent
paragraphs $X$ is an arbitrary complex manifold and $L$ a line bundle
on it.

Let $\pi: \Lbullet\to X$ be the $\C^*$ bundle associated to $L$, as
follows: on $\Lbullet$ over $U_i$ we have a coordinate $z_i\in\C^*$
and $z_i=g_{ij}z_j$.  In other words $\Lbullet={\rm
  Spec}_X(\bigoplus_{k\in\Z}L^{\otimes k})$, or equivalently,
$\Lbullet$ is the total space of $L^*$ with zero section removed.  On
$\Lbullet$ we have a natural action of $\C^*$, that is: $\C^*\times
\Lbullet\ni(t,z)\mapsto t\cdot z\in \Lbullet$ and related
$\C^*$-invariant non-vanishing vector field $(z_i/\partial z_i)$
trivializes the relative tangent bundle $T_{\Lbullet/X}$ which is the
kernel of the tangential map $T\pi: T_{\Lbullet}\to T_X$. Its dual
$\Omega^1_{\Lbullet/X}$ is thus trivialized by a $\C^*$-invariant form
$\mu=dz_i/z_i$. So we have a short exact sequence
$$
0\to \pi^*(\Omega^1_X) \to \Omega^1_{\Lbullet} \to
\Omega^1_{{\Lbullet}/X} \cong \O_{\Lbullet} \to 0
$$
which we can push down to $X$ and use the fact $\pi_*(\O_\Lbullet)
= \bigoplus L^{\otimes k}$ to get
$$
0 \to \bigoplus_{k\in\Z}\Omega^1_X(k) \to \pi_*(\Omega^1_{\Lbullet})
\cong \bigoplus_{k\in\Z} (\Omega^1_{\Lbullet})_k \to \bigoplus_{k\in\Z}
L^{\otimes k} \to 0
$$
where $\pi_*(\Omega^1_{\Lbullet})\cong\bigoplus
(\Omega^1_{\Lbullet})_k$ is the splitting into the weight spaces of
the $\C^*$-action. We denote $(\Omega^1_{\Lbullet})_0$ by ${\mathcal
  L}$ and thus find an exact sequence, called the Atiyah extension of
$L$:
$$
0 \to \Omega^1_X \to {\mathcal L} \to \O_X \to 0
$$

\begin{fact}[Atiyah]
  The Chern class $c_1(L)\in H^1(X,\Omega^1_X) =
  \Ext^1_X(\O_X,\Omega^1_X)$ is the extension class of the above
  sequence.
\end{fact}

\begin{rem}
  The sheaf ${\mathcal L}(1)$ may be identified as the sheaf of first
  jets of sections of $L$, see e.g.~\cite[1.6.3]{BS95}, but we do not
  use this fact.
\end{rem}

\begin{fact}
If $f: Z\to X$ is a morphism then $c_1(f^*(L))$ is the image of
$c_1(L)$ in the composition
$$
c_1(L)\in \Ext_X^1(\O_X,\Omega^1_X)\to
\Ext^1_Z(\O_Z,f^*(\Omega^1_X))\to \Ext^1_Z(\O_Z,\Omega^1_Z)\ni
c_1(f^*(L))
$$
with the first arrow coming from pull-back and the second one from
the derivative $df: f^*(\Omega^1_X)\to\Omega^1_Z$. In particular, if
$c_1(f^*(L))$ is non-zero then $f^*({\mathcal L})$ is coming from a
non-trivial extension in $\Ext^1_Z(\O_Z,f^*(\Omega^1_X))$.
\end{fact}

\subsubsection{Projectivized Vector Bundles}
Let us consider a rank $r+1$ vector bundle ${\mathcal E}$ over a
smooth manifold $Y$. The projectivization of ${\mathcal E}$, denoted
by $\P({\mathcal E})$, is defined as the relative projective spectrum
${Proj}_X({Sym}({\mathcal E}))$ which geometrically means the variety
of lines in the dual bundle ${\mathcal E}^*$. The variety
$\P({\mathcal E})$ admits a projection $p: \P({\mathcal E})\to Y$,
with fibers being projective spaces $\P_r$, and a line bundle
$\O_{\P({\mathcal E})}(1)$ whose restriction to each fiber is
$\O_{\P_r}(1)$. Furthermore, we have $p_*(\O_{\P({\mathcal E})}(1)) =
{\mathcal E}$. The cokernel of the differential map $dp:
p^*(\Omega^1_Y) \to \Omega^1_{\P({\mathcal E})}$ is
$\Omega^1_{\P({\mathcal E})/Y}$, the relative cotangent bundle.  

The functoriality of the first Chern class implies the following
lemma, whose proof we omit referring the reader to \cite[II,
8.13]{Ha77} for a definition of the Euler sequence.

\begin{lem}\label{jar:2.1.2}
  The image of the extension class $c_1(\O_{\P({\mathcal E})}(1))\in
  \Ext^1_{\P({\mathcal E})}(\O_{\P({\mathcal E})},\Omega^1_{\P({\mathcal
      E})})$ under the induced map 
  $$
  \Ext^1_{\P({\mathcal E})}(\O_{\P({\mathcal
      E})},\Omega^1_{\P({\mathcal E})}) \to \Ext^1_{\P({\mathcal
      E})}(\O_{\P({\mathcal E})},\Omega^1_{\P({\mathcal E})/Y})
  $$
  is associated with the (dual) relative Euler sequence
  $$
  0\to \Omega^1_{\P({\mathcal E})/Y} \to p^*({\mathcal
    E})\otimes\O_{\P({\mathcal E})}(-1) \to \O_{\P({\mathcal E})}\to
  0.
  $$
\end{lem}

\subsection{Symplectic structure on $\Lbullet$}
For an introduction to these matters, see \cite[prop.~2.4]{Leb95}. Now
we assume that $X$ is a contact manifold with contact form $\theta$.
Then on ${\Lbullet}$ we define a weight 1 homogeneous form $\omega\in
H^0(X,(\Omega^2_{\Lbullet})_1)\cong H^0(X, \wedge^2{\mathcal L}(1))$
by setting:
$$
\omega=d\theta+\theta\wedge\mu= d\theta_i+\theta_i\wedge
(dz_i/z_i)
$$
It is easy to verify that this gives a global 2-form on ${\bf
  L}^\bullet$ and moreover that $\omega^{n+1} =
(n+1)\cdot\theta\wedge\mu\wedge (d\theta)^n$ so that the form $\omega$
is symplectic. Thus we obtain an isomorphism $\omega: {\mathcal
  L}^*\to {\mathcal L}(1)$ which we denote by the same name.

\subsection{Splitting type of $T_X$ on rational curves}
It is a well-known fact attributed to Grothendieck that any vector
bundle $E$ on $\P_1$ decomposes into a sum of line bundles, $E\cong
\bigoplus \O(e_i)$ where the (non-increasing) sequence of integers
$(e_i)$ is called the splitting type of $E$.

The main point in the proof of theorem~\ref{mainthm} is an analysis of
the splitting type of $T_X$, restricted to minimal rational
curves. For this, we introduce the following

\begin{notation}
If $E$ is a vector bundle on $\P_1$, let $\rk+(E)$ denote the number
of positive entries in the splitting type, that is,
$\rk+(E)=\#\{e_i>0\}$.
\end{notation}

We have the following elementary lemma:

\begin{lem}
  Suppose that a bundle $E$ is in an exact sequence of bundles on
  $\P_1$: 
  $$
  0\to E\to \bigoplus_{i=1}^r \O(a_i)\oplus\bigoplus_{j=1}^s
  \O(-b_j) \to \O\to 0
  $$
  with $r\geq 0$, $s>0$ and all $a_i\geq 0$, $b_j>0$. If the
  sequence does not split, then $\rk+(E^*)=s-1$.
\end{lem}
\begin{proof}
  Since the sequence does not split, the map $H^0(E) \to
  H^0(\bigoplus_{i=1}^r \O(a_i))$ is an isomorphism and thus the
  number of non-negative values in the splitting type of $E$ is $r$.
  Indeed, let $E^+\hookrightarrow E$ be the sub-bundle associated to
  the non-negative values in the splitting type: then
  $H^0(E^+)=H^0(E)=H^0(\bigoplus\O(a_i))$ and the embedding
  $E^+\hookrightarrow \bigoplus\O(a_i)$ is an isomorphism.  Therefore
  the number of negative values in the splitting type of $E$ is $s-1$.
\end{proof}

\begin{prop}\label{jar:2.2.2}
  Let $f:\P_1\to C\subset X$ be a normalization of a rational curve in
  the contact manifold $X$. If $\deg(f^*(L))=1$, then $\rk+(f^*(T_X))
  = n$
\end{prop}
\begin{proof}
  The symplectic form $\omega$ on $\Lbullet$ gives an isomorphism
  $f^*({\mathcal L}^*)\cong f^*({\mathcal L})\otimes\O(1)$, see
  \cite[prop.~2.4]{Leb95}. It follows that
  $$
  f^*({\mathcal L})= \bigoplus_{i=1}^{n+1}
  (\O(a_i)\oplus\O(-a_i-1))
  $$
  where $a_i$ are non-negative. On the other hand, by the preceding
  discussion on the Atiyah extension of $L$, it follows that the
  pullback of the sequence defining ${\mathcal L}$ to $\P_1$ does not
  split. Therefore we can apply the previous lemma and obtain that
  $\rk+(f^*(T_X))=n$.
\end{proof}

\subsection{Extremal rational curves on contact manifolds}
In the present section we apply proposition~\ref{jar:2.2.2} to study
the locus of minimal rational curves passing through a given point of
$X$. We use Koll\'ar's book \cite{Kol96} as our main reference.

Let us consider the scheme $\Hom(\P_1,X)$ parameterizing morphisms
from $\P_1$ to $X$. By $\Hom(\P_1,X;0\mapsto x)$ we denote the scheme
parameterizing morphisms sending $0\in \P_1$ to $x\in X$.  Let
$F:\P_1\times \Hom(\P_1,X)\to X$ be the evaluation morphism. The
restriction of $F$ to any family of morphisms $V\subset \Hom(\P_1,X)$
will be denoted by $F_V$. We write $\locus(V)$ for the image of the
evaluation $F_V(\P_1\times V)$. By abuse, we will say that $V$ is a
family of rational curves and $\locus(V)$ is their locus. If
$x\in\locus(V)$ is any point, then let $\locus(V,x)$ denote
$\locus(V\cap \Hom(\P_1,X; 0\mapsto x))$, that is, the locus of curves
from $V$ which contain $x$.

Following Koll\'ar, we say that the family $V$ is {\em unsplit} if
curves from the family can not be deformed to a 1-cycle consisting of
more than one (counting possible multiplicities) curves ---for a
precise definition we refer to \cite[IV, 2.1]{Kol96}. We note that in
our case $V$ is unsplit if $f(\P_1)$ is an extremal rational curve in
the sense of Mori (see the definition below).

\begin{prop}\label{jar:2.3.1}
  Let $X$ be a projective contact manifold of dimension $2n+1$ and let
  $V$ be an irreducible component of $\Hom(\P_1,X)$.  If $V$ is an
  unsplit family of rational curves and for $f\in V$ we have
  $\deg(f^*(L))=1$, then $\locus(V)=X$ and moreover $\dim
  \locus(V,x)=n$ for any $x\in X$.
\end{prop}

\begin{proof}
  Let $f:\P_1\to X$ be a morphism such that $[f]\in V$ and $f(0)=x$.
  By \cite[II 3.10]{Kol96} and proposition~\ref{jar:2.2.2} the
  tangential map $TF$ satisfies
  $$
  \rank \left[TF_{(p,[f])}\right]=\rk+ f^*(T_X)=n
  $$
  for any $p\in \P_1-\{0\}$. Therefore, by \cite[III, 10.6]{Ha77},
  for {\em any} $x\in \locus(V)$ we have $\dim \locus(V,x)\leq n$. On
  the other hand, by Koll\'ar's version of the Ionescu-Wi\'sniewski
  estimate \cite[IV, 2.6.1]{Kol96} for any $x\in\locus(V)$ we have
  $$
  \dim X+\deg(f^*(-K_X))=3n+2\leq \dim\locus(V)+\dim \locus(V,x)+1
  $$
  Note that the quoted inequality \cite[IV, 2.6.1]{Kol96}) is
  stated for a {\em general} point $x$ in the locus of a {\em
    generically} unsplit family $V$. However we claim that it remains
  true for {\em any} $x\in \locus(V)$ if the family is unsplit.
  The proof is the same as the one of \cite[prop.~IV, 2.5]{Kol96}
  ---we point out that there is an obvious misprint reversing the
  inequality in this proposition. Now comparing the inequalities we
  get:
  $$
  3n+2\leq \dim\locus(V)+\dim\locus(V,x)+1 \leq \dim X +
  n+1=3n+2
  $$
  so that $\locus(V)=X$ and $\dim\locus(V,x)=n$ for any
  $x\in X$.
\end{proof}

\subsection{Contractions of contact manifolds. Proof
  of theorem 1.1}
\label{sect:contr_of_rays}

For the convenience of the reader we recall basic facts of Mori theory
before commencing the proof of theorem~\ref{mainthm}. For a complete
treatment of this subject we refer to e.g.~\cite{KM98}.

\subsubsection{Mori Theory}
Let $X$ be a complex projective manifold such that $K_X$ is not nef.
This means that there exists a curve $C$ such that $K_X\cdot C<0$.
Then, by the Kawamata-Shokurov base-point-free theorem, $X$ admits a
Mori contraction. That is, there exists a normal projective variety
$Y$ and a surjective morphism $\phi: X \to Y$ (which is not an
isomorphism) such that $\phi_*(\O_X)=\O_Y$ (i.e.~its fibers are
connected) and $-K_X$ is $\phi$-ample.

The contraction $\phi$ is called {\em elementary} if all curves
contracted by $\phi$ to points are numerically proportional, or
equivalently, if $b_2(X)=b_2(Y)+1$. It is a basic fact that any
contraction can be factored via an elementary one.

If $\phi$ is an elementary contraction, then we define its {\em
  length} $l(\phi):=\min\{-K_X\cdot C\}$, where $C$ is among rational
curves contracted by $\phi$ (a rational curve $C$ for which the
minimum is achieved we call an extremal rational curve). According to
Mori's cone theorem, the number $l(\phi)$ is defined (i.e.~there
exists a rational curve contracted by $\phi$) and $\dim X+1\geq
l(\phi)\geq 1$.

If $X$ is a projective contact manifold then, since $-K_X=(n+1)L$, the
bundle $L$ is $\phi$-ample for any Mori contraction $\phi$ of $X$.
Moreover, if $\phi$ is elementary then $l(\phi)$ is either $2n+2$ or
$n+1$, and the latter occurs if there exists a rational curve $C$
contracted by $\phi$ such that $L\cdot C=1$. If $l(\phi)=2n+2=\dim
X+1$ then, by a result of Ionescu, see e.g.~\cite[IV, 2.6.1]{Kol96},
the contraction $\phi$ is onto a point and therefore $X$ is Fano with
$b_2=1$.

\subsubsection{Proof of theorem~\ref{mainthm}}

In view of the theory outlined in the previous paragraph, in order to
prove theorem~\ref{mainthm}, it remains to show that an elementary
contraction $\phi: X\to Y$ is always isomorphic to $\P(T_Y)\to Y$.

As a first step, we show that in our setup Mori contractions cannot be
birational.

\begin{lem}\label{jar:2.4.1}
  Let $\phi: X\to Y$ be a Mori contraction of a projective contact
  manifold $X$. Then $\dim Y<\dim X$.
\end{lem}
\begin{proof}
  Since any Mori contraction can be factored through an elementary one
  we may assume that $\phi$ is elementary and moreover, by the
  previous considerations, that $l(\phi)=n+1$, because otherwise $Y$
  is a point. Let $C\subset X$ be a rational curve contracted by
  $\phi$ such that $-K_X\cdot C=(n+1)\cdot L\cdot C=n+1$. Let $f:
  \P_1\to C$ be a normalization of $C$ and consider an irreducible
  component $V$ of $\Hom(\P_1,X)$ which contains $[f]$. By
  proposition~\ref{jar:2.3.1}, $\locus(V)=X$ so that the lemma
  follows.
\end{proof}

\begin{prop}\label{jar:2.4.2}
  Let $X$ be as above and $\phi :X\to Y$ be a surjective morphism to a
  variety where $0<\dim Y<\dim X$. If $X_\eta$ is a generic fiber and
  $X_\eta$ is Fano, then $X_\eta \cong \mathbb P_n$ and $X_\eta$ is an
  integral manifold of $F$, i.e.~$T_{X_\eta}\subset F$.
\end{prop}
\begin{proof}
We may assume that $X_\eta$ and $\phi (X_\eta)$ are smooth and
that $\phi $ has maximal rank at every point of $X_\eta$.

As a first step we show that $X_\eta \cong \mathbb P_n$. In order to
do this, construct a sheaf-morphism $\beta :L|_{X_\eta}\to
T_{X_\eta}$.  Take $Hom(.,L)$ of the sequence defining the contact
structure and obtain 
$$
\begin{CD}
  0 @>>> \O_X @>>>\Omega^1_X\otimes L @>>> F^*\otimes L @>>> 0
\end{CD}.
$$
Contraction with the contact form $\omega \in H^0(X,(F\otimes
F)^*\otimes L)$ yields a morphism $\iota \omega :F\to F^*\otimes L$.
Since $\omega $ is non-degenerate, this must be an isomorphism $F\cong
F^*\otimes L$. Thus, one obtains a map $\alpha :\Omega ^1_X\otimes
L\to T_X$ as follows:
\begin{equation}
  \begin{CD}
    0 @>>> \O_X @>>>\Omega^1_X\otimes L @>>>
    F^*\otimes L @>>> 0 \\ @. @. @VV{\alpha}V @VV{\cong}V \\ 0 @<<< L @<<<
    T_X @<<< F @<<< 0
  \end{CD}\label{Fiber:comdiag1}
\end{equation}
Now consider the canonically defined map $\phi ^*(\Omega ^1_Y)\otimes
L\to \Omega ^1_X\otimes L$.  Restrict this map to $X_\eta$ and note
that $\phi^*(\Omega _Y^1)|_{X_\eta}\cong \oplus ^{\dim
  Y}\O_{X_\eta}$.  This yields a sequence of morphisms as follows:
$$
\xymatrix { & & T_{X_\eta}\ar[d] \\
  (\phi^*(\Omega^1_Y)\otimes L)|_{X_\eta} \cong ({L^{\oplus \dim
      Y})|_{X_\eta}} \ar[r] \ar@/^/@{-->}[urr]^{\beta} &
  (\Omega^1_X\otimes L)|_{X_\eta}\ar[r]^{\alpha|_{X_\eta}} &
  T_X|_{X_\eta}} 
$$
The induced map $\beta $ exists because the normal bundle
$N_{X_\eta}$ is trivial, but $L|_{X_\eta}$ is ample so that every
morphism $L\to N_{X_\eta}$ is necessarily trivial.

Now claim that $\beta$ is not the trivial map. Since $\phi$ has
maximal rank at some point of $X_\eta$, the map $(L^{\oplus
  m})|_{X_\eta}\to (\Omega _X^1\otimes L)|_{F}$ is not identically
zero. On the other hand, diagram \ref{Fiber:comdiag1} implies directly
that $Ker(\alpha )\cong \O _X$. Again using that there is no map
$L|_{X_\eta}\to \O _{X_\eta}$, we see that $\beta$ is non-trivial
indeed.

The existence of $\beta $ is equivalent to
$h^0(X_\eta,T_{X_\eta}\otimes L^*|_{X_\eta})\geq 1$. In this situation
a theorem of J.~Wahl \cite{Wah83} applies, showing that $X_\eta\cong
\mathbb P_k$ for some $k\in \mathbb N$. Use the adjunction formula to
see that $k=n$.

In order to see that $X_\eta$ is $F$-integral, let $x\in X_\eta$ be
any point and note that the tangent space $T_{X_\eta}|_{\{x\}}$ is
spanned by the tangent space of the lines in $X_\eta$ passing through
$x$. Thus, in order to show that $X_\eta$ is integral, it is
sufficient to show that all lines in $X_\eta$ are integral manifolds.
If $C \subset X_\eta$ is a line, using the adjunction formula, we
obtain $-K_X.C=n+1$, i.e.~$L.C=1$. If $\iota : C\to X$ is the
embedding, note that, since $T_C\cong \O_C(2)$, the induced map
$T_C\to L$ is necessarily trivial. This shows that $T_{X_\eta}\subset
F|_{X_\eta}$.
\end{proof}

The preceding proposition can also be proved using the Atiyah extension class of $L$ and the
symplectic form $\omega,$ c.f. (2.2). We may now finish with the proof of theorem~\ref{mainthm}.

\begin{thm}\label{jar:2.4.4}
  Let $\phi: X\to Y$ be a Mori contraction of a projective contact
  manifold $X$ of dimension $2n+1$ onto a positive dimensional variety
  $Y$.  Then $Y$ is smooth of dimension $n+1$, moreover $X\cong
  \P(T_Y)$ and $L\cong \O_{\P(T_Y)}(1)$.
\end{thm}
\begin{proof}
  First we prove that $Y$ is smooth and $X=\P(\phi_*(L))$.  This will
  follow from a result by Fujita \cite[2.12]{Fuji85}, if only we prove
  that the morphism $\phi$ is equidimensional (i.e.~it has all fiber
  of dimension $n$).
  
  After twisting $L$ by a pull-back of an ample line bundle from $Y$
  we may assume that the result, call it $L'$, is ample. Let ${\rm
    Chow}_{n,1}(X)$ be a Chow variety of $n$-dimensional cycles on $X$
  of degree 1 with respect to $L'$, for the definition see
  e.g.~\cite[I, 3]{Kol96}.  Let $X_g$ be a general fiber of $\phi$, by
  proposition~\ref{jar:2.4.2} $(X_g,L'_{| X_g})\cong (\P_n,\O(1))$. We
  consider an irreducible component $W\subset {\rm Chow}_{n,1}(X)$
  which contains $[X_g]$. We note that $W$ is projective of dimension
  $n+1$ and over $W$ there exists a universal family of cycles $\pi:
  U\to W$ with the evaluation map $e: U\to X$.  The map $e$ is
  birational and either it is an isomorphism or it has a
  positive-dimensional fiber. In the former case, however, $(\phi:
  X\to Y)\cong (\pi: U\to W)$ and thus $\phi$ is equidimensional and
  we are done by Fujita's result. Thus, to arrive to a contradiction,
  we assume that $e$ has a positive dimensional fiber. We note that,
  if $x_0\in X$ is such that $\dim(e^{-1}(x_0))>0$, then by the
  property of functor $Chow$, we have
  $\dim\left(e(\pi^{-1}(\pi(e^{-1}(x_0))))\right) > n$, because
  otherwise all cycles in $\pi(e^{-1}(x_0))$ would be mapped to one
  $n$-dimensional cycle.
  
  Let us consider an irreducible component $V_U$ of $\Hom(\P_1, U)$
  such that $f(\P_1)$ is a line in $X_g=\pi^{-1}([X_g])\cong\P_n$ and
  $f$ is isomorphism onto its image. The family $V_U$ is unsplit and
  it maps, via the natural map $\tilde e: \Hom(\P_1,U)\to{\rm
    Hom}(\P_1,X)$, into an unsplit component $V$ of $\Hom(\P_1,X)$
  containing curves whose degree with respect to $L$ is 1. Now, we
  note that $\locus(V_U,u)=\pi^{-1}(\pi(u))\cong \P_n$ for a general
  $u\in U$, and because $V_u$ is unsplit the equality $\locus(V_U,u) =
  \pi^{-1} (\pi(u))$ holds for any $u\in U$. Therefore $\locus(V,x_0)
  \supset e(\pi^{-1}(\pi(e^{-1}(x_0))))$ and since the dimension of
  the latter set is bigger than $n$ we arrive to a contradiction to
  proposition~\ref{jar:2.3.1}. This concludes the proof of the first
  part of theorem~\ref{jar:2.4.4}.
  
  Thus $Y$ is smooth and we set ${\mathcal E}:=\phi_*(L)$ so that 
  $(X,L)=(\P({\mathcal E}),\O_{\P({\mathcal E})}(1))$. We shall prove
  that ${\mathcal E}\cong T_Y$. By lemma~\ref{jar:2.1.2} we have a map of
  extensions
  $$
  \xymatrix{ 
    0 \ar[r]& {\Omega^1_X}\ar[r]\ar[d] &
    {\mathcal L}\ar[r]\ar[d] & {\O_X} \ar[r]\ar@{=}[d] &0\\
    0\ar[r]& {\Omega^1_{X/Y}}\ar[r]&
    {\phi^*{\mathcal E}(-1)}\ar[r]& {\O_X}\ar[r]& 0
    }
  $$
  and comparing it with the exact sequence defining $\Omega^1_{X/Y}$
  we get an exact sequence
  $$
  0\to\phi^*(\Omega^1_Y) \to {\mathcal L} \to \phi^*{\mathcal E}(-1)
  \to 0
  $$
  Now we combine the embedding $\phi^*(\Omega^1_Y) \to {\mathcal L}$
  with the symplectic isomorphism $\omega: {\mathcal L}\cong{\mathcal
    L}^*(-1)$.  By considering the restriction to fibers of $\phi$ we
  see that the induced map, coming from the twisted dual of the
  extension defining ${\mathcal L}$
  $$
  \begin{CD}
    \phi^*(\Omega^1_Y) @>>> {\mathcal L} @>{\omega}>> {\mathcal L}^*(-1)
    @>>> \phi^*((\Omega^1_Y)^*)(-1)
  \end{CD}
  $$
  is zero so we have an embedding of $\phi^*(\Omega^1_Y)$ into the
  kernel of ${\mathcal L}^*(-1) \to \phi^*((\Omega^1_Y)^*)(-1)$ which is
  $\phi^*({\mathcal E}^*)$. The resulting map is of maximal rank so
  that we get $\phi^*(\Omega^1_Y) \cong \phi^*({\mathcal E}^*)$ which
  implies ${\mathcal E}=T_Y$ and concludes the proof of
  theorem~\ref{jar:2.4.4}.
\end{proof}

\subsection{Contact structures on $\P(T_Y)$}

We now turn to contact structures on projectivized tangent bundles
which we shall describe in more detail. At first, we consider contact
varieties with more than one elementary contraction. It was known
before that if $X$ is Fano, then either $b_2(X)=1$ or
$X=\P(T_{\P_{n+1}})$. This was proved in \cite[cor.~4.2]{LBS94} using
Wi\'sniewski's classification of Fano manifolds. In fact, a stronger
version of that result holds:

\begin{prop}\label{prop:two_contrs}
  A projective contact manifold $X$ admits one extremal ray
  contraction at most, unless $X\cong \P(T_{\P_{n+1}})$.
\end{prop}
\begin{proof}
  If $X$ has more than one extremal ray then, by the Mori cone
  theorem, the cone of curves of $X$ has also a Mori face of dimension
  $\geq 2$. This means that there exists a Mori contraction $\psi:
  X\to Y$ with $b_2(Y)\leq b_2(X)-2$ which can be factored through at
  least two different elementary contractions $\phi_i$ as follows:
  $$
  \xymatrix {
    X \ar[r]^{\phi_1}\ar[d]_{\phi_2}\ar[dr]^{\psi} & {Y_1} \ar[d]\\
    {Y_2} \ar[r] & Y
    }
  $$
  We claim that in this case $Y$ is just a point so we are in the
  case of a Fano contact manifold with $b_2\geq 2$ which, by the
  mentioned above result, is $\P(T_{\P_{n+1}})$.  Indeed, otherwise
  $\dim Y>0$ and we can apply proposition~\ref{jar:2.4.2} to obtain
  $\psi^{-1}(y)\cong\P_n$, for a general $y\in Y$. But by
  lemma~\ref{jar:2.4.1} $\psi^{-1}(y)$ must contain curves contracted
  by both $\phi_1$ and $\phi_2$, a contradiction, since $b_2(\P_n)=1$
  and thus all curves on $\psi^{-1}(y)\cong\P_n$ are numerically
  proportional.
\end{proof}

If $X=\P(T_{\P_{n+1}})$, then the contact structure on $X$ is unique, see
e.g. \cite[cor.~3.2]{Leb95}. Note also that $\Aut(\P(T_{\P_{n+1}})) /
\Aut(\P_{n+1}) = \{1\}$. We will see that this is
not a coincidence.

Let $\phi: X=\P(T_Y)\to Y$ be the projectivization of the tangent bundle on a compact complex manifold $Y$ of dimension $n+1$, together with the relative ample sheaf $L=\O_{\P(T_Y)}(1)$. We have an exact sequence of bundles of twisted differentials $$
0 \to \phi^*(\Omega^1_Y)(1) \to \Omega^1_X(1) \to \Omega^1_{X/Y}(1) \to 0 $$
We claim that the resulting map
$$
H^0(X,\phi^*(\Omega^1_Y)(1))\to H^0(X,\Omega^1_X(1)) $$
is an isomorphism. This follows from observing that $H^0(X,\Omega^1_{X/Y}(1))$ is the kernel of the evaluation map $$
H^0(X,\phi^*(T_Y))=H^0(X,\phi^*(\phi_*(\O_{\P(T_Y)}(1)))\to H^0(X,\O_{\P(T_Y)}(1))
$$
appearing in the cohomology of twisted relative Euler sequence $$
0 \to \Omega^1_{X/Y}(1) \to \phi^*(T_Y) \to \O_{\P(T_Y)}(1) \to 0 $$
On the other hand we have naturally, by push-forward, $$
H^0(X,\phi^*(\Omega^1_Y)(1)) = H^0(Y,\phi_*(\phi^*(\Omega^1_Y)(1))) = H^0(Y, \Omega^1_Y\otimes T_Y)=H^0(Y, \End(\Omega^1_Y)), $$
so sections of $\phi^*\Omega^1_Y(1)$ can be identified naturally with endomorphisms of $\Omega^1_Y$.

Let us reveal the local nature of the above identification. In an analytic neighborhood $U$ of a point $y\in Y$ we choose 
local (analytic !) coordinates $(y^\alpha)$ and, accordingly, we choose local vector fields $\zeta^\beta$, where $0\leq \alpha, \beta\leq n$, 
so that $\zeta^\alpha(y^\alpha)=1$ and $\zeta^\beta(y^\alpha)=0$ if $\alpha\ne\beta$. We note that $\zeta^\beta$ are also 
homogeneous coordinates in $\phi^{-1}(U)$ and, equivalently, generating sections of $\O_{\P(T_Y)}(1)_{|\phi^{-1}(U)}$. 
If $h=(h^{\beta\alpha})\in \End(\Omega^1_Y)(U)\cong Mat((n+1)\times (n+1), \O(U))$ is an endomorphism of $\Omega^1_Y$ over $U$, 
represented in bases related to the choice of coordinates $(y^\alpha)$ as a $(n+1)\times (n+1)$ matrix then setting
$$
\theta_h=\sum_{0\leq \alpha, \beta\leq n} \zeta^\beta h^{\beta\alpha} dy^\alpha
$$
we get a local presentation of $h$, related to the choice coordinates, a well as a well defined 1-form in $\Omega^1_X(1)_{|\phi^{-1}(U)}$.

\begin{prop}\label{jar:2.5.1}
The above defined identification
$$
H^0(Y,\End(\Omega^1_Y))\ni h\mapsto \theta_h\in H^0(X,\Omega^1_X(1))
$$
provides a bijection between automorphisms of $\Omega^1_Y$ and contact forms on $X=\P(T_Y)$.
\end{prop}

\begin{proof}
It is enough to prove that $\theta_h\wedge (d\theta_h)^n\in H^0(X,K_X\otimes L^{n+1})$ does not vanish anywhere if and only if $h$ is an 
automorphism. Since, however, $K_X\otimes L^{n+1}=\O_X$ it is enough to verify the non-vanishing at one point. This will be done using 
the local description provided above. First, we compute $$
d\theta_h=\sum_{0\leq\alpha, \beta\leq n} \left[d\zeta^\beta\wedge(h^{\beta \alpha}d y^{\alpha}) + 
\zeta^{\beta}\cdot (dh^{\beta \alpha}\wedge dy^{\alpha})\right] $$
and we see that terms of type $\zeta^{\beta}\cdot (dh^{\beta 
\alpha}\wedge dy^{\alpha})$ can be ignored in the computation of $\theta_h\wedge(d\theta_h)^n$ because they do not 
contain differentials with respect to $\zeta^\beta$. On the other hand we note that on the projective space $\P_n$ with 
homogeneous coordinates $[\zeta^0,\dots,\zeta^n]$ we can identify the form $\sum(-1)^k\cdot\zeta^k\cdot
d\zeta^0\wedge\cdots\hat{{\phantom{.}^k\phantom{.}}} \dots\wedge d\zeta^n$ with the unit section of $\Omega_{\P_n}^n(n+1)\cong \O_{\P_n}$. 
Moreover, if write $d\tilde y^\beta=\sum_\alpha h^{\beta\alpha} dy^\alpha$ then 
$d\tilde y^0\wedge\cdots\wedge d\tilde y^n=\det(h)\cdot dy^0\wedge\cdots\wedge dy^n$. So, taking all the above into account and remembering 
that the wedge product of degree 2 forms commutes we get
\begin{align*}
&\theta_h\wedge(d\theta_h)^n \\
&= n!\cdot \sum_{k=0\dots n} (\zeta^k\cdot d\tilde y^k)\wedge[ (d\zeta^0\wedge d\tilde
y^0)\wedge\cdots\hat{{\phantom{.}^k\phantom{.}}} \cdots\wedge(d\zeta^n\wedge d\tilde y^n)] \\ &=\pm n!\sum_{k=0\dots n} 
(\zeta^k\cdot d\tilde y^k)\wedge (d\zeta^0\wedge\cdots\hat{\phantom{.}^k\phantom{.}}\cdots\wedge d\zeta^n)\wedge (d\tilde
y^0\wedge\cdots\hat{\phantom{.}^k\phantom{.}}\cdots\wedge d\tilde y^n) \\
&=\pm n!\sum_{k=0\dots n} (-1)^k\cdot \zeta^k\cdot (d\zeta^0\wedge\cdots\hat{\phantom{.}^k\phantom{.}}\cdots\wedge d\zeta^n)\wedge
(d\tilde y^0\wedge\cdots\wedge d\tilde y^n) \\ &=\pm n! \cdot \det(h) \cdot \left[\sum_k (-1)^k\cdot \zeta^k\cdot 
(d\zeta^0\wedge\cdots\hat{\phantom{.}^k\phantom{.}}\cdots\wedge d\zeta^n)\right] \wedge(dy^0\wedge\cdots\wedge dy^n) \end{align*}
and the last expression is non-zero if and only if $\det(h)\ne 0$. This concludes the proof of proposition~\ref{jar:2.5.1}. \end{proof}

\section{Contact manifolds where $K_X$ is nef}
\label{KxNEF_CHAP}

In this section we turn to the case where $K_X$ is nef. It is shown in
\cite{Dru98} that a projective contact manifold always has Kodaira
dimension $\kappa (X)=-\infty$. The abundance conjecture predicts that
this is incompatible with $K_X$ nef. However, since the conjecture is
known only in dimension $\leq 3$ and completely open in higher
dimensions, we have to consider this possibility here.

As a partial result, we show in theorem~\ref{2Cprop} that $K_X$ is not
nef if $X$ has more than one contact structure, and in
proposition~\ref{Knef:stability} that $K_X$ is not nef if a certain
stability property of $T_X$ holds.

\begin{prop}\label{2Cprop}
  Let $X$ be a projective manifold admitting at least two contact
  structures. Then the canonical bundle $K_X$ is not nef and thus
  either $X=\P (T_Y)$ or $X$ is Fano and $b_2(X)=1$.
\end{prop}

\begin{proof}
  Let
  \begin{equation}
    \label{Contact_Sequence}
    0\to F\to T_X\to L\to 0
  \end{equation}
  and
  $$
  0\to F'\to T_X\to L'\to 0
  $$
  be two different contact structures on $X$. Since
  $L^{n+1}=(L')^{n+1}=-K_X$, the line bundle $L-L'$ is torsion and,
  after passing to a finite \'etale cover, we may assume that $L=L'$.
  
  The map $F\to L'$ yields a non-zero element $v\in H^0(F^*\otimes
  L')=H^0(F^*\otimes L)\cong H^0(F)\subset H^0(T_X)$, i.e.~we obtain a
  vector field. If $v$ has zeroes, then $X$ is uniruled and we see
  that $K_X$ cannot be nef. Thus, we assume that $v$ has no zeroes. In
  this situation \cite[thm.~3.13]{Lib78} asserts that ---after another
  \'etale cover, if necessary--- $X\cong T\times Y$, where $T$ is a
  torus and $v\in \pi ^*_{1}(H^0(T_{T}))\subset H^0(F)$. Here
  $\pi_1:X\to T$ and $\pi_2 : X \to Y$ are the natural projections. It
  follows from the adjunction formula that the restriction of $L$ to
  $\pi_2$-fibers is torsion.  Thus, again taking covers, we may assume
  that $L=\pi^*_2 (L_Y)$ with $L_Y\in \Pic (Y)$.  Let $U\subset Y$
  denote an open set such that $L_Y|_U=\O_U$.  Let
  $$
  \theta \in H^0(\pi^{-1}_2(U),\Omega^1_X\otimes L) \cong
  H^0(\pi^{-1}_2(U),\Omega_X^1) \cong H^0(\pi^*_1(\Omega_T^1)) \oplus
  H^0(\pi^*_2(\Omega_Y^1))
  $$
  be the contact form associated with
  sequence~\ref{Contact_Sequence} and $\theta = \theta_T \oplus
  \theta_Y$ be the direct sum decomposition. By assumption, the
  following expression is not zero:
  $$
  \theta \wedge (d\theta )^{\wedge n}=(\theta_T+\theta_Y)\wedge
  (d\theta_T+d\theta_Y)^{\wedge n}=\theta_T\wedge
  (d\theta_Y)^{\wedge n}.
  $$
  This, however, is absurd: first, it follows from
  $(d\theta_Y)^{\wedge n}\ne 0$ that $\dim Y=2n$, $\dim T=1$. It then
  follows from $\theta (v)=0$ that $\theta_T=0$. A contradiction.
\end{proof}

As a first step towards the proof of proposition~\ref{Knef:stability},
the succeeding result asserts that in our situation $K^2_X$ is zero.
For this we do not actually need that $X$ has a contact structure and
consider a more general situation:

\begin{prop}\label{Knef:KX_square}
  Let $X$ be a projective manifold. Let $L^*\subset \Omega _X^1$ be a
  locally free subsheaf of rank 1 with $\alpha L^*\equiv K_X$ for some
  positive rational number $\alpha $. If $K_X$ is additionally nef,
  then $K^2_X=0$.
\end{prop}
The notation $L^*$ might seem to be slightly awkward at first. We use
it to be consistent with the notation introduced in
section~\ref{morichap:jar}.

\begin{proof}
  Consider the smooth surface $S$ cut out by general hyperplane
  sections
  $$
  S:=H_1\cap \ldots \cap H_{\dim X-2}.
  $$
  Since $L^*$ is nef, there is no morphism $L^*|_S\to N^*_{S|X}$,
  and we obtain an injection $L^*\to \Omega_S^1$. By Bogomolov's
  well-known theorem (see \cite[thm.~on p.~501]{Bog79}), $L^*|_S$
  cannot be big. Thus $(L^*)^2.S = 0$ which translates into
  $$
  K^2_X.H_1\ldots H_{\dim X-2}=0.
  $$
  We claim that if $D$ is any nef divisor with $D^2.H_1\ldots
  H_{\dim X-2}=0.$ for generic ample divisors $H_i$, then $D^2\equiv
  0$.  If $\dim X = 3$, the claim is obvious. If $\dim X\geq 5$, then
  we can apply the Lefschetz theorem and reduce to the case where
  $\dim X=4$ by taking suitable hyperplane sections.
  
  If $\dim X=4$ and $H$ is any ample divisor, then $D+\epsilon H$ is
  an ample $\mathbb R$-divisor for all $\epsilon >0$. Consequently,
  since $(D+\epsilon H)^2.H_1\in \overline{NE(X)}$, we have
  $D^2.H_1\in \overline{NE(X)}$, and we conclude that $D^2.H_1\equiv
  0$. Now apply the Hodge index theorem in the form of
  \cite[App.~A, 5.2]{Ha77} with $Y := D^2$ and $H := H_1$. Assuming
  that $D\not \equiv 0$ that theorem yields that $D^4>0$. But
  $$
  D^4 = \lim_{\epsilon \to 0} D^2.\underbrace{(D+\epsilon
    H_1).(D+\epsilon H_2)}_{\txt{\scriptsize two gen. ample divisors}}
  = 0,
  $$
  a contradiction.
\end{proof}
The proof the preceding proposition gives rise to the following

\begin{problem}\label{Knef:Problem}
  Let $S$ be a smooth projective surface and $L^*\subset \Omega^1_X$ a
  numerically effective locally free subsheaf of rank one with
  $(L^*)^2=0$ and $L^*\not \equiv 0$. Can anything be said about the
  structure of $S$? Observe that if $L^*\subset \Omega_X^1$ is a
  subbundle away from a finite set (e.g.~if $X$ is a contact
  manifold), then $L^*|_S\subset \Omega_S^1$ is a subbundle away from
  a finite set as well.
\end{problem}
\begin{cor}
  In the setting of proposition \ref{Knef:KX_square}, the tangent
  bundle $T_X$ is semistable with respect to $(K_X,H^{\dim X-2})$,
  where $H$ is any ample line bundle. This means
  $$
  \frac1{r}c_1(\mc S)\cdot K_X\cdot H^{\dim X-2}\leq \frac1{\dim
    X}c_1(X)\cdot K_X\cdot H^{\dim X-2}
  $$
  for all coherent subsheaves $\mc S\subset T_X$ of any rank $r>0$.
\end{cor}
\begin{proof}
  Since $K^2_X\equiv 0$, this follows from \cite{Eno87}.
\end{proof}

Unfortunately it is not always true that $T_X$ is also semistable with respect to $(K_X+\epsilon
H,H^{m-2})$ for some ample $H$ and some small number $\epsilon >0$. Counterexamples
are e.g. provided by products of elliptic curves with curves of genus $g \geq 2.$  
On the other hand, a weak consequence of
this semistability property would already imply the assertion that
$K_X$ is not nef:

\begin{thm}\label{Knef:stability}
  In the setting of proposition \ref{Knef:KX_square}, if there exists
  an ample bundle $H\in \Pic(X)$ and a number $\epsilon >0$ such that
  $L^*$ does not destabilize $\Omega^1_X$ with respect to
  $(K_X+\epsilon H,H^{\dim X-2})$ and if $0 < \alpha < \dim X$, then
  either $K_X$ is not nef or $K_X\equiv 0$.
\end{thm}
\begin{proof}
  The assertion on the destabilization can be expressed as
  $$
  c_1(L^*)\cdot (K_X+\epsilon H)\cdot H^{\dim X-2}\leq \frac1{\dim
    X}K_X\cdot (K_X+\epsilon H)\cdot H^{\dim X-2}.
  $$
  Using $K^2_X=0$ and $K_X=\alpha L^*$, the inequality becomes
  \begin{equation}\label{Knef:stability_ineq}
    \frac1{\alpha }K_X\cdot H^{\dim X-1}\leq \frac1{\dim X}K_X\cdot 
    H^{\dim X-1}
  \end{equation}
  If we now assume that $K_X$ is nef, then
  inequality~(\ref{Knef:stability_ineq}) implies $K_X\cdot H^{\dim
    X-1}=0$, which is equivalent to $K_X\equiv 0$.
\end{proof}

Therefore it remains to consider the case that $L^*$ destabilizes $T_X$
for all polarizations $(K_X + \epsilon H,H^{\dim X-2}).$ One might
hope to derive some geometric consequences from this unstability.

Remark that all considerations of the present section apply in
particular to contact manifolds where $\dim X=2n+1$ and $\alpha =n+1$.
Thus, if the unstable case could be handled,
then the canonical bundle $K_X$ of a projective contact manifold $X$
would never be nef. This is because the assumption $K_X\equiv 0$
contradicts a result of Ye \cite[lem. 1]{Ye94}.

\section{Manifolds with nef subsheaves in $\Omega _X^1$}
\label{nef_in_Omega_CHAP}

The setup of proposition~\ref{Knef:KX_square} seems to be of
independent interest. The aim of the present section is to give a
description of these varieties.

\begin{assumption}\label{subbdle:setting}
  Let $X$ be a projective manifold and $L^*\subset \Omega_X^1$ be a
  locally free nef subsheaf of rank one. Assume that there is a
  positive rational number $\alpha$ such that $\alpha L^*=K_X$.
\end{assumption}
Recall the well-known result of Bogomolov \cite{Bog79} which implies
that in this setting $\kappa (X)\leq 1$.

As a first result we obtain:

\begin{prop}
  If the Kodaira dimension $\kappa (X)\geq 0$, then $K_X\equiv 0$ or
  $\alpha \geq 1$.
\end{prop}
\begin{proof}
  We argue by absurdity: assume that $K_X \not \equiv 0$ and $\alpha <
  1$. The inclusion $L^*\to \Omega_X^1$ gives a non-zero element
  $$
  \theta \in H^0(\Omega_X^1\otimes L)=H^0\left(\wedge^{\dim X-1}T_X\otimes
    \left( 1-\frac1{\alpha}\right) K_X\right)
  $$
  Suppose that $\alpha <1$. Then $(1-\frac1{\alpha })<0$ and thus
  we can find positive integers $m$, $p$ such that
  $$
  H^0(\underbrace{S^{m}(\wedge^{\dim X-1}T_X)\otimes \O (-pK_X)}_{=:E})
  \ne 0.
  $$
  In order to derive a contradiction, recall the result of Miyaoka
  (\cite[cor.~8.6]{Miy87}; note that $X$ cannot be uniruled) that
  $\Omega_X^1|_C$ is nef for a sufficiently general curve $C\subset X$
  cut out by general hyperplane sections of large degree. Remark that
  $E^*|_C$ is nef, too. We may assume without loss of generality that
  $p$ is big enough so that $pK_X$ has a section with zeroes, and so
  does $E|_C$. See \cite[prop.~1.2]{CP91} for a list of basic
  properties of nef vector bundles which shows that this is
  impossible.
\end{proof}

\subsection{The case where \protect$\kappa (X)=1\protect$}

In this case $K_X$ is semi-ample and we give a description of
$L^*$ in terms of the Kodaira-Iitaka map.

\begin{thm}\label{Knef:Iitaka}
  Under the assumptions~\ref{subbdle:setting}, if $\kappa (X)=1$, then
  $K_X$ is semi-ample, i.e.~some multiple is generated by global
  sections.  Let $f:X\to C$ be the Iitaka fibration and $B$ denote the
  divisor part of the zeroes of the natural map $f^*(\Omega_C^1)\to
  \Omega_X^1$.  Then there exists an effective divisor $D\in \Div(X)$
  such that $L^*=f^*(\Omega^1_C)\otimes \O _X(B-D)$.
  
  Furthermore, if $p$ is chosen such that $pK_X\in f^*(\Pic(C))$ and
  such that $p\alpha$ is an integer, then
  $$
  p(\alpha -1)K_C=f_*(pK_{X|C})+f_*(p\alpha (D-B)).
  $$
\end{thm}
\begin{proof}
  Since $K_X\not \equiv 0$ and $K^2_X=0$ by proposition
  \ref{Knef:KX_square}, the numerical dimension $\nu (X)=1$. In this
  setting \cite[thm.~1.1]{Kaw85} proves that $K_X$ is semi-ample,
  i.e.~that a sufficiently high multiple of $K_X$ is globally
  generated. By \cite[lem. 12.7]{Bog79}, the canonical morphism
  $$
  L^*\to \Omega _{X|C}^1
  $$
  is generically 0. Let $B$ denote the divisor part of the zeroes
  of $f^*(\Omega_{C}^1)\to \Omega_X^1$.  Then we obtain an exact
  sequence
  $$
  0\to f^*(\Omega _{C}^1)\otimes \O _X(B)\to \Omega _X^1\to
  \tilde{\Omega }_{X|C}^1\to 0
  $$
  where the cokernel $\tilde{\Omega}_{X|C}^1$ is torsion free.
  Since $\Omega ^1_{X|C}=\tilde{\Omega }^1_{X|C}$ away from a closed
  subvariety, the induced map
  $$
  L^*\to \tilde{\Omega }_{X|C}^1
  $$
  vanishes everywhere. Hence there exists an effective divisor $D$
  such that
  $$
  L^*\subset f^*(\Omega _{C}^1)\otimes \O _X(B)\quad \textrm{i.e.}
  \quad L^*=f^*(\Omega _{C}^1)\otimes \O _X(B-D).
  $$
  The last equation is obvious.
\end{proof}

\subsection{The case where \protect$\kappa (X)=0\protect$}

We now investigate the more subtle case $\kappa(X)=0$. We pose the
following

\begin{conjecture}\label{subbdle:our_conjecture}
  Under the assumptions~\ref{subbdle:setting}, if $\kappa (X)=0$,
  then $K_X\equiv 0$. Hence (see \cite{Bea83}) there exists a finite
  \'etale cover $\gamma :\tilde X\to X$ such that $\gamma^*(L^*) =
  \O_{\tilde X}$ and $\tilde X=A\times Y$ where $A$ is Abelian and $Y$
  is simply connected.
\end{conjecture}
We will prove the conjecture in a number of cases, in particular if
$\dim X\leq 4$ or if the well-known Conjecture K holds (see
\cite[sect.~10]{Mor87} for a detailed discussion). Recall that the
Albanese map of a projective manifold $X$ with $\kappa(X)=0$ is
surjective and has connected fibers \cite{Kaw81}.

\begin{conjecture}[Ueno's Conjecture K]\label{subbdle:conjecture_k}
  If $X$ is a nonsingular projective variety with $\kappa (X)=0$, then
  the Albanese map is birational to an \'etale fiber bundle over
  $\Alb(X)$ which is trivialized by an \'etale base change.
\end{conjecture}
It is known that Conjecture K holds if $q(X)\geq \dim X-2$.

The proof requires two technical lemmata. The first is a
characterization of pull-back divisors. It appears implicitly in
\cite[p.~571]{Kaw85}.

\begin{lem}[Kawamata's pull-back lemma]\label{subbdle:pull-back-lemma}
  Let $f:X\to Y$ be an equidimensional surjective projective morphism
  with connected fibers between quasi-projective $\mathbb Q$-factorial
  varieties and let $D\in \mathbb Q\Div(X)$ be an $f$-nef $\mathbb
  Q$-divisor such that $f(\Supp(D))\not =Y$. Then there exists a
  number $m\in \mathbb N^{+}$ and a divisor $A\in \Div(Y)$ such that
  $mD=f^*(A)$.
\end{lem}
The next lemma is very similar in nature, and also the proof is very
kind. We include it here for lack of an adequate reference.

\begin{lem}\label{subbdle:lambda_negative}
  Let $f:X\to Y$ be a surjective projective morphism with connected
  fibers between quasi-projective $\mathbb Q$-factorial varieties and
  let $D\in \mathbb Q\Div(X)$ be an $f$-nef $\mathbb Q$-divisor which
  is of the form
  $$
  D=f^*(A)+\sum \lambda_iE_i
  $$
  where $A\in \mathbb Q\Div(Y)$ and $\codim_Y f(E_i)\geq 2$.
  Then $\lambda_i \leq 0$ for all $i$.
\end{lem}
\begin{proof}
  Suppose that this was not the case. Choose $j$ such that $\lambda_j
  > 0$ and such that $f(E_j)$ is of maximal dimension. As the lemma is
  formulated for quasi-projective varieties, in order to derive a
  contradiction, we may assume without loss of generality that $f(E_i)
  = f(E_j)$ for all $i$.
  
  \subsubsection*{Claim:} There exists a divisor $B\in \mathbb 
  Q\Div(Y)$ and a number $m\in \mathbb N^+$ such that the $\mathbb
  Q$-divisor $M:=mD-f^*(B)$ satisfies the following conditions
  \begin{enumerate}
  \item $-M$ is effective
  \item there exists a component $F_i$ of $F:=f^{-1}f(E_j)$ which is
    not contained in $M$
  \item $F\cap \Supp(M)\not =\emptyset$
  \end{enumerate}
  \subsubsection*{Application of the claim:} If the claim is true, 
  then we can always find a curve $C\subset F$ such that $C\not
  \subset \Supp(M)$, $C\cap \Supp(M)\not =\emptyset $ and $f(C)=(*)$.
  Since $-M$ is effective, we have $C.M<0$, contradicting the assumed
  $f$-nefness.
  
  \subsubsection*{Proof of the claim:} choose a very ample effective
  Cartier-divisor $H\in \Div(Y)$ such that
  \begin{itemize}
  \item the effective part of $A$ is contained in $H$ and
  \item $f(E_i)\subset \Supp(H)$. 
  \end{itemize}
  If $M:=D$ satisfies the requirements of the claim already, stop
  here. Otherwise, set
  $$
  D':=(\underbrace{\textrm{mult. of $E_j$ in
      $f^{*}(H)$}}_{\textrm{positive}}).D-\underbrace{(\textrm{mult.
      of $E_j$ in D})}_{\textrm{positive
      }}.\underbrace{f^{*}(H)}_{\textrm{cont. $E_i$ with pos. mult.}}
  $$
  by choice of $H$ and by choice of the coefficients, $D'$
  satisfies conditions (2) and (3) of the claim. If $-D'$ is
  effective, then we are finished already.
  
  If $-D'$ is not effective, note that the positive part of $D'$ is
  supported on $F$ only. Repeat the above procedure with a new number
  $j$.  Note that after finitely many steps the divisor must become
  anti-effective. This finishes the proof.
\end{proof}
We will now show that Conjecture K implies our
conjecture~\ref{subbdle:our_conjecture}.

\begin{thm}\label{subbdle:conj_k_implies_K_triv}
If conjecture K holds, then conjecture \ref{subbdle:our_conjecture}
holds as well.
\end{thm}
In order to clarify the structure of the proof, we single out the case
where $H^0(X,L^*)\not =0$. More precisely, we consider the
weakened

\begin{assumption}\label{subbdle:modified_setting}
  Let $X$ be a projective manifold and $L^*\subset \Omega_X^1$ be a
  locally free nef subsheaf of rank one. Assume that $\kappa (X)=0$
  and that there is a positive rational number $\alpha $ such that
  $\alpha L^*\subset K_X$ is a subsheaf. Assume furthermore that
  $H^0(X,L^*)\not =0$.
\end{assumption}

\begin{lem}\label{subbdle:sect_and_K}
  Under the weakened assumptions~\ref{subbdle:modified_setting}, if
  Conjecture K holds for $X$, then
  conjecture~\ref{subbdle:our_conjecture} holds as well,
  i.e.~$L^*\equiv 0$.
\end{lem}
\begin{proof}
  In inequality $H^0(X,L^*)\not =0$ directly implies that $q(X)>0$,
  i.e.~that the Albanese map $f:X\to A = \Alb(X)$ is not trivial.
  
  In this setting, Conjecture K yields a diagram as follows: 
  $$
  \xymatrix { & {X_{dom}} \ar[dl]_{\pi} \ar[dr]^{\pi '}\\ X \ar[d]_f
    \ar@{-->}[rr] & & X' \ar[d]_{\txt{\scriptsize \'etale \\ 
        \scriptsize fiber\\ \scriptsize bundle}}^{f'}& & {F\times
      \tilde A}\ar[d]^{\pi_2} \ar[ll]_{\Gamma}^{\txt{\scriptsize
        \'etale, finite}} \ar[r]^{\pi_1} & F \\ A \ar@{=}[rr] & & A &
    & {\tilde{A}} \ar[ll]_{\gamma}^{\txt{\scriptsize \'etale, finite}}
    }
  $$
  where $\pi$ and $\pi'$ are birational morphisms. Let $E$ be the
  zero-set of a section of $L^*$; i.e.~$L^*=\O (E)$. We will show that
  $E=0$.
  
  \subsubsection*{Claim 1:} The divisor $E$ is supported on 
  $f$-fibers, i.e.~$f(\Supp(E))\not =A$.
  
  \subsubsection*{Proof of claim 1:} the inclusion $L^*\subset 
  \Omega _X^1$ yields an element
  $$
  \theta \in H^0(\Omega _X^1(-E)).
  $$
  Since $\omega = f^*(\eta)$, we conclude that
  $$
  f(E)\subset \Sing(f),
  $$
  where $\Sing(f)$ is the set of all $y\in A$ such that $f^{-1}(y)$
  is singular. In particular, $E$ does not meet the general fiber of
  $f$. This shows claim 1.
  
  As a next step in the proof of lemma~\ref{subbdle:sect_and_K}, we
  set
  $$
  A^0:=\{a\in A|\dim f^{-1}(a)=\dim X-\dim A\}\quad \textrm{ and
    }\quad X^0:=f^{-1}(A^0)
  $$
  and show
  
  \subsubsection*{Claim 2:} the divisor $E$ does not intersect $X^0$,
  i.e.~$\Supp(E)\cap X^0=\emptyset$.
  
  \subsubsection*{Proof of claim 2:} by Kawamata's pull-back 
  lemma~\ref{subbdle:pull-back-lemma} there exists a number $m\in
  \mathbb N^{+}$ such that \( mE|_{X^0}\in f^*(D|_{A^0})$ for some
  $D\in \Div(X)$. This implies
  $$
  \overline{f^{-1}(D\cap A^0)}\subset \Supp(E).
  $$
  Now consider $X_{dom}$. Since $X$ is smooth, we can find a
  $\mathbb Q$-divisor for the canonical bundle $\omega _{X_{dom}}$ of
  the form
  $$
  K_{X_{dom}}=\pi ^*(\underbrace{K_X}_{\alpha
    E+(\textrm{effective})})+\sum \lambda_i B_i
  $$
  where the $B_{i}$ are $\pi $-exceptional divisors and
  $\lambda_{i}\in \mathbb Q^{+}$.  In particular, we have
  $$
  \overline{(\pi \circ f)^{-1}(D\cap X^0)}\subset
  \Supp(\pi^*(E))\subset \Supp(K_{X_{dom}}).
  $$
  Since $\pi \circ f=\pi '\circ f'$ and $K_{X'} :=
  (\pi')_*(K_{X_{dom}})$ is a $\mathbb Q$-divisor for the
  anticanonical bundle $\omega _{X'}$, we find
  $$
  \overline{(f')^{-1}(D\cap X^0)}\subset \Supp(K_{X'}).
  $$
  We derive a contradiction from the last inequality by noting that
  $$
  (\pi_2\circ \gamma )^{-1}(D)\subset \Supp(\underbrace{K_{F\times
      \tilde{A}}}_{=\Gamma ^*(K_{X'})}).
  $$
  On the other hand, it follows from the adjunction formula that
  every effective $\mathbb Q$-divisor for the canonical bundle
  $\omega_{F\times \tilde{A}}$ must be contained in $\pi^*_1(\mathbb
  Q\Div(F))$. This shows claim 2.
  
  \subsubsection*{Application of claims 1 and 2:} we know that $E$ 
  is an effective and nef divisor satisfying $\codim_A(f(\Supp(E)))
  \geq 2$. By lemma~\ref{subbdle:lambda_negative} this is possible if
  and only if $E=0$.
\end{proof}
The preceding lemma enables us to finish the proof of
theorem~\ref{subbdle:conj_k_implies_K_triv} by reducing to the
modified weakened assumptions~\ref{subbdle:modified_setting}. Since we
wish to re-use the same argumentation later, we formulate a technical
reduction lemma:

\begin{lem}[Reduction Lemma]\label{subbdle:reduction_lemma} 
  If $L^*\equiv 0$ holds true for all varieties of a given dimension
  $d$ satisfying the weakened
  assumptions~\ref{subbdle:modified_setting}, then
  conjecture~\ref{subbdle:our_conjecture} holds for all varieties of
  dimension $d$.
\end{lem}
\begin{proof}
  Let $X$ be a variety as in conjecture \ref{subbdle:our_conjecture}
  and assume that $X$ is of dimension $d$. If $H^0(X,L^*)\not =0$,
  then we can stop here.

  Otherwise, we find a minimal positive integer $m$ and an effective
  divisor $E$ such that $(L^*)^{m}=\O _X(E)$. We build a sequence of
  morphisms
  $$
  \begin{CD}
    \tilde Y@>{\sigma}>{\txt{\scriptsize desing.}}>Y
    @>{\gamma}>{\txt{\scriptsize cyclic cover}}>\hat X
    @>{\pi}>\txt{\scriptsize desing. of $E$}> X 
  \end{CD}
  $$
  as follows: Let $\pi :\hat X\to X$ be a sequence of blow-ups with
  smooth centers such that $\Supp(\pi ^*(E))$ has only normal
  crossings.  Set $\hat{E}=\pi^*(E)$ and $\hat{L}=\pi^*(L)$ and let
  $f:Y\to \hat X$ be the cyclic covering (followed by normalization)
  associated with the section
  $$
  s\in H^0(X,(L^*)^{m})
  $$
  which defines $\hat{E}$; see e.g.~\cite[3.5]{EV92} for a more
  detailed description of this construction. By \cite[3.15]{EV92}, $Y$
  is irreducible, \'etale over $\hat X\setminus \Supp(\hat{E})$ and
  smooth over $\hat X\setminus \Sing(\Supp(\hat{E}))$. Moreover
  $$
  H^0(Y,\gamma ^*(\hat{L}^*))\not =0.
  $$
  Let $\sigma :\tilde Y\to Y$ be a desingularization. We will
  show:
  
  \subsubsection*{Claim: $\kappa (\tilde Y)=0$.}
  
  \subsubsection*{Application of the claim:} if the claim holds true, we
  conclude as follows: let $\tilde{L}=(\sigma\circ\gamma)^*(\hat{L})$.
  Then $H^0(\tilde X,\tilde{L}^*)\not =0$, and furthermore
  $\tilde{L}^*\subset \Omega_{\tilde Y}^1$ by virtue of the canonical
  morphism
  $$
  (\underbrace{\sigma \circ \gamma \circ \pi}_{=: \Gamma}
  )^*(\Omega_X^1) \to \Omega_{\tilde Y}^1.
  $$
  Similarly, since $X$ is smooth, $\Gamma^* K_X \subset K_{\tilde
    Y}$.  Thus, lemma~\ref{subbdle:sect_and_K} applies to $\tilde Y$,
  showing that $\Gamma^* (L^*)\equiv 0$ so that $L^*$ was numerically
  trivial in the first place. This shows the lemma and thus finishes
  the proof of theorem~\ref{subbdle:conj_k_implies_K_triv}.
  
  \subsubsection*{Proof of the claim:} Since $\Gamma$ is an
  \'etale cover away from $E$, we can find a divisor for the canonical
  bundle $\omega_{\tilde Y}$ which is of the form
  $$
  K_{\tilde Y}=\Gamma ^*(\alpha E)+D
  $$
  where $D$ is effective and $\Gamma (D)\subset E$. This already
  implies that there is a number $k\in \mathbb N^+$ such that
  $(k\Gamma^*(E)-K_{\tilde Y})$ is effective. Thus $\kappa (K_{\tilde
    Y})\leq \kappa (\Gamma^*(E))=0$, and the claim is shown. 
\end{proof}
Thus ends the proof of theorem~\ref{subbdle:conj_k_implies_K_triv}.

Finally we show conjecture \ref{subbdle:our_conjecture} in the case
where $\dim X\leq 4$.

\begin{thm}\label{subbdle:dimX_small}
  Under the assumptions~\ref{subbdle:setting}, if $\kappa (X)=0$ and
  $\dim X\leq 4$, then $K_X\equiv 0$.
\end{thm}
\begin{proof}
  Using the reduction lemma~\ref{subbdle:reduction_lemma}, it is
  sufficient to consider the weakened
  setting~\ref{subbdle:modified_setting}: let $X$ be as
  in~\ref{subbdle:modified_setting}. At first we argue exactly as in
  the proof of lemma \ref{subbdle:sect_and_K}: the Albanese map
  $f:X\to A=\Alb(X)$ is surjective and has connected fibers and we have
  $f(\Supp(E))\not =A$, where $E$ is the zero-divisor of the section
  in $L^*$. Furthermore we have $\kappa (X,E)=0$.
  
  Recall that Conjecture K holds if $q(X)\geq \dim X-2$; see
  \cite[p.~316]{Mor87}. By lemma~\ref{subbdle:sect_and_K} we are
  finished in these cases. It remains to consider the case where
  $q(X)=1$.
  
  Because $f$ is equidimensional in this case, Kawamata's pull-back
  lemma \ref{subbdle:pull-back-lemma} applies: there is a number $m\in
  \mathbb N^{+}$ such that $mE=f^*(D)$. It follows immediately that
  $\kappa (A,D)=0$.  Since $A$ is a torus, this is possible if and
  only if $D=0$.
\end{proof}

\begin{rem}
  Actually, the proof of lemma~\ref{subbdle:reduction_lemma} and
  theorem~\ref{subbdle:dimX_small} show that
  conjecture~\ref{subbdle:our_conjecture} holds if $q(X)\geq \dim X-2$
  or if $H^0(X,L^*)\not =0$ and $q(X)=1$. For the first statement,
  note that $q(X)$ increases when passing to a cover.
\end{rem}


\begin{thebibliography}{KPS99}

\bibitem[Ati57]{Ati57}
M.~Atiyah.
\newblock {Complex analytic connections in fibre bundles}.
\newblock {\em Trans. Am. Math. Soc.}, 85:181--207, 1957.

\bibitem[Bea83]{Bea83}
A.~Beauville.
\newblock {Vari\'{e}t\'{e}s K\"{a}hl\'{e}riennes dont la premi\`{e}re classe
  de Chern est nulle}.
\newblock {\em J. Diff. Geom.}, 18:755--782, 1983.

\bibitem[Bea98]{Bea98}
A.~Beauville.
\newblock Fano contact manifolds and nilpotent orbits.
\newblock {\em Comm. Math. Helv.}, 73(4):566--583, 1998.

\bibitem[Bea99]{Bea99}
A.~Beauville.
\newblock Riemannian holonomy and Algebraic Geometry.
\newblock Duke/alg-geom Preprint 9902110, 1999.

\bibitem[Bog79]{Bog79}
F.~Bogomolov.
\newblock Holomorphic tensors and vector bundles on projective varieties.
\newblock {\em Math. USSR Izv.}, 13:499--555, 1979.

\bibitem[BS95]{BS95}
M.~Beltrametti and A.~Sommese.
\newblock {\em The Adjunction Theory of Complex Projective Varieties}.
\newblock de Gruyter, 1995.

\bibitem[CP91]{CP91}
F.~Campana and T.~Peternell.
\newblock Projective manifolds whose tangent bundles are 
          numerically effective.
\newblock {\em Math. Ann.}, 289:169--187, 1991.

\bibitem[Dru98]{Dru98}
S.~Druel.
\newblock Contact structures on algebraic 5-dimensional manifolds.
\newblock C.R. Acad. Sci.Paris, 327:365-368, 1998.

\bibitem[Eno87]{Eno87}
I.~Enoki.
\newblock Stability and negativity for tangent bundles of minimal K\"ahler
  spaces.
\newblock In T.~Sunada, editor, {\em Geometry and Analysis on Manifolds.
  Proceedings, Katata-Kyoto}, volume 1339 of {\em Lecture Notes in
  Mathematics}, pages 118--127. Springer, 1987.

\bibitem[EV92]{EV92}
H.~Esnault and E.~Viehweg.
\newblock {\em Lectures on Vanishing Theorems}.
\newblock Birkh\"{a}user Verlag, 1992.

\bibitem[Fuj85]{Fuji85}
T.~Fujita.
\newblock {On Polarized Manifolds Whose Adjoint Bundles Are Not Semipositive}.
\newblock In T.~Oda, editor, {\em Advanced Studies in Pure Mathematics 10},
  pages 167--178. North-Holland publishing company, 1985.

\bibitem[Har77]{Ha77}
R.~Hartshorne.
\newblock {\em Algebraic Geometry}, volume~52 of {\em Graduate Texts in
  Mathematics}.
\newblock Springer, 1977.

\bibitem[Kaw81]{Kaw81}
Y.~Kawamata.
\newblock {Characterization of Abelian varieties}.
\newblock {\em Compositio Math.}, 43:253--276, 1981.

\bibitem[Kaw85]{Kaw85}
Y.~Kawamata.
\newblock Pluricanonical systems on minimal algebraic varieties.
\newblock {\em Inv. Math.}, 79:567--588, 1985.

\bibitem[Kol96]{Kol96}
J.~Koll\'ar.
\newblock {\em Rational Curves on Algebraic Varieties}.
\newblock volume 32 of {\em Ergebnisse der Mathematik und ihrer Grenzgebiete}. Springer, 1996.

\bibitem[KM98]{KM98}
J.~Koll\'ar and S.~Mori.
\newblock Birational Geometry of algebraic varieties.
\newblock  volume 134 of {\em Cambridge Tracts in Math}. Cambridge
University Press, 1998.

\bibitem[LeB95]{Leb95}
C.~LeBrun.
\newblock Fano manifolds, contact structures and quaternionic geometry.
\newblock {\em Int. Journ. of Math.}, 6(3):419--437, 1995.

\bibitem[LS94]{LBS94}
C.~LeBrun and S.~Salamon.
\newblock Strong rigidity of positive quaternion-K\"{a}hler manifolds.
\newblock {\em Inv. Math.}, 118(1):109--132, 1994.

\bibitem[Lie78]{Lib78}
D.~Liebermann.
\newblock Compactness of the Chow scheme: Applications to automorphisms and
  deformations of K\"ahler manifolds.
\newblock In Francois Norguet, editor, {\em Fonctions de Plusieurs Variables
  Complexes III}, number 670 in Lecture Notes in Mathematics, pages 140--186.
  Springer, 1978.

\bibitem[Miy87]{Miy87}
Y.~Miyaoka.
\newblock Deformation of a morphism along a foliation.
\newblock In S.~Bloch, editor, {\em Algebraic Geometry}, volume~46(1) of {\em
  Proceedings of Symposia in Pure Mathematics}, pages 245--268, Providence,
  Rhode Island, 1987. American Mathematical Society.

\bibitem[Mor87]{Mor87}
S.~Mori.
\newblock Classification of higher-dimensional varieties.
\newblock In {\em Proceedings of Symposia in Pure Mathematics}, volume~46, of
{\em Proceedings of Symposia in Pure Mathematics}, pages 269--331,
Providence, Rhode Island, 1987. American Mathematical Society.


\bibitem[Wah83]{Wah83}
J.~Wahl.
\newblock {A cohomological characterization of $\mathbb P_n$}.
\newblock {\em Inv. Math.}, 72:315--322, 1983.

\bibitem[Ye94]{Ye94}
Y.-G. Ye.
\newblock A note on complex projective threefolds admitting holomorphic 
  contact structures.
\newblock {\em Inv. Math.}, 115:311--314, 1994.

\end{thebibliography}
\end{document}